\documentclass[11pt]{article}
\usepackage{amsmath, amssymb, amsthm, verbatim,enumerate,bbm,color} %DIF >
\numberwithin{equation}{section}
\usepackage{float}

\usepackage{mathrsfs}

\usepackage{tikz}
\usetikzlibrary{calc}
\usetikzlibrary{math}
\usetikzlibrary{decorations,decorations.markings,decorations.pathreplacing,decorations.pathmorphing}
\usetikzlibrary{fadings}
\usetikzlibrary{fit}
\usetikzlibrary{matrix}
\usetikzlibrary{patterns}
\usetikzlibrary{positioning}
\usetikzlibrary{shadows}
\usetikzlibrary{shapes,shapes.geometric}
\usetikzlibrary{trees}
\usetikzlibrary{overlay-beamer-styles}
\usepackage{tikzsymbols}

\usepackage[colorlinks,citecolor=blue,bookmarks=true]{hyperref}

\date{\today}
\allowdisplaybreaks

\theoremstyle{plain}
\newtheorem{theorem}{Theorem}[section]

\newtheorem{proposition}[theorem]{Proposition}

\newtheorem{problem}[theorem]{Problem}

\makeatletter
\def\moverlay{\mathpalette\mov@rlay}
\def\mov@rlay#1#2{\leavevmode\vtop{%
		\baselineskip\z@skip \lineskiplimit-\maxdimen
		\ialign{\hfil$\m@th#1##$\hfil\cr#2\crcr}}}
\newcommand{\charfusion}[3][\mathord]{
	#1{\ifx#1\mathop\vphantom{#2}\fi
		\mathpalette\mov@rlay{#2\cr#3}
	}
	\ifx#1\mathop\expandafter\displaylimits\fi}
\makeatother

\makeatletter
\renewenvironment{proof}[1][\proofname]
{\par\pushQED{\qed}
	\normalfont\topsep6\p@\@plus6\p@\relax\trivlist
	\item[\hskip\labelsep\bfseries#1\@addpunct{.}]
	\ignorespaces}
{\popQED\endtrivlist\@endpefalse}

\makeatother

\RequirePackage{color}\definecolor{RED}{rgb}{1,0,0}\definecolor{BLUE}{rgb}{0,0,1} %DIF PREAMBLE
\newcommand{\poly}{\text{poly}}

\usepackage{hyperref}

\begin{document}

\title{A Generalization of Varnavides's Theorem}

\author{Asaf Shapira\thanks{School of Mathematics, Tel Aviv University, Tel Aviv 69978, Israel. Email: asafico@tau.ac.il. Supported in part
by ERC Consolidator Grant 863438 and NSF-BSF Grant 20196.}}

\date{}
\maketitle

\begin{abstract}

A linear equation $E$ is said to be {\em sparse} if there is $c>0$ so that every subset of $[n]$ of size $n^{1-c}$ contains a solution of $E$ in distinct integers. The problem of characterizing the sparse equations, first raised by Ruzsa in the 90's, is one of the most important open problems in additive combinatorics.
We say that $E$ in $k$ variables is {\em abundant} if every subset of $[n]$ of size $\varepsilon n$ contains at least $\poly(\varepsilon)\cdot n^{k-1}$ solutions of $E$. It is clear that every abundant $E$ is sparse, and
Gir\~{a}o, Hurley, Illingworth and Michel asked if the converse implication also holds. In this note we show that this is the case for every $E$ in $4$ variables. We further discuss a generalization of this problem which applies to all linear equations.

%Fix a linear equation $E$ of the form $\sum^k_{i=1}a_ix_i=0$ with $\sum_ia_i=0$.
%One of the most well studied problems in additive combinatorics
%asks to determine $R_E(n)$, which denotes the size of the largest subset of $[n]$ not containing a solution of $E$ in distinct integers.
%Of particular importance is an old problem of Ruzsa which asks to determine for which $E$ there is $c=c(E)$ so that $R_E(n) \leq n^{1-c}$.
%Such equations are said to be {\em easy}. It is also
%natural to let $S_E(\varepsilon)$ be such that every subset of $[n]$ of size $\varepsilon n$ contains at least $n^{k-1}/S_{E}(\varepsilon)$ %solutions of $E$. Clearly $S_{E}(\varepsilon) \leq (1/\varepsilon)^C$ implies that $E$ is easy, and
%Gir\~{a}o et al. recently asked if the converse implication also holds. In this note we show that this is the case if $E$ is an equation in %$k=4$ variables.
\end{abstract}

\section{Introduction}\label{sec:intro}

Tur\'an-type questions are some of the most well studied problems in combinatorics.
They typically ask how ``dense'' should
an object be in order to guarantee that it contains a certain small substructure. In the setting of graphs,
this question asks how many edges an $n$-vertex graph should contain in order to force the appearance of some fixed graph $H$.
For example, a central open problem in this area asks, given a bipartite graph $H$, to determine the smallest $T=T_{H}(\varepsilon)$ so that for every $n\geq T$ every
$n$-vertex graph with $\varepsilon {n \choose 2}$ edges contains a copy of $H$
(see \cite{Bukh} for recent progress). A closely related question which also attracted a lot of attention, is the {\em supersaturation}
problem, introduced by Erd\H{o}s and Simonovits \cite{ES} in the 80's. In the setting of Tur\'an's problem for bipartite
$H$, the supersaturation question asks to determine the largest $T^*_{H}(\varepsilon)$ so that every $n$-vertex graph with $\varepsilon {n \choose 2}$ edges contains
at least $(T^*_{H}(\varepsilon)-o_n(1)) \cdot n^h$ labelled copies of $H$, where $h=|V(H)|$ and $o_n(1)$ denotes a quantity tending to $0$ as $n \rightarrow \infty$. One of the central conjectures in this area, due to Sidorenko, suggests
that $T^*_{H}(\varepsilon) = \varepsilon^{m}$, where $m=|E(H)|$ (see \cite{CLS} for recent progress).

We now describe two problems in additive number theory, which are analogous to the graph problems described above.
We say that a homogenous linear equation $\sum^k_{i=1}a_ix_i=0$ is {\em invariant} if $\sum_ia_i=0$. All equations we consider here will be invariant and homogenous.
Given a fixed linear equation $E$, the Tur\'an problem for $E$ asks to determine
the smallest $R=R_{E}(\varepsilon)$ so that for every $n\geq R$, every $S\subseteq [n]:=\{1,\ldots,n\}$ of size $\varepsilon n$ contains a solution to $E$ in distinct integers.
For example, when $E$ is the equation $a+b=2c$ we get the Erd\H{o}s--Turán--Roth problem on sets avoiding $3$-term arithmetic progressions (see \cite{KM} for recent progress).
Continuing the analogy with the previous paragraph, we can now ask to determine the largest $R^*_{E}(\varepsilon)$
so that every $S \subseteq [n]$ of size $\varepsilon n$ contains at least $(R^*_{E}(\varepsilon)-o_n(1)) \cdot n^{k-1}$ solutions to $E$, where $k$ is the number of variables in $E$.
We now turn to discuss two aspects which make the arithmetic problems more challenging than the graph problems.

Let us say that $E$ is {\em sparse} if there is $C=C(H)$ so that\footnote{It is easy to see that this definition is equivalent to the one we used in the abstract.} $R_E(\varepsilon) \leq \varepsilon^{-C}$. The first aspect which makes the arithmetic landscape more varied is that while in the case of graphs it is well known (and easy) that for every bipartite $H$ we have $T_{H}(\varepsilon)=\poly(1/\varepsilon)$, this is no longer the case in the arithmetic setting. Indeed, while Sidon's equation $a+b=c+d$ is sparse, a well known construction of Behrend \cite{Behrend} shows that
$a+b=2c$ is not sparse\footnote{More precisely, it shows that in this case $R_{E}(\varepsilon) \geq (1/\varepsilon)^{c\log 1/\varepsilon}$. Here and throughout this note, all logarithms are base $2$.}.
%$R_{E}(\varepsilon) \geq (1/\varepsilon)^{c\log 1/\varepsilon}$.
The problem of determining which equations $E$ are sparse is a wide open problem due to Ruzsa, see Section $9$ in \cite{R}.

Our main goal in this paper is to study another aspect which differentiates the arithmetic and graph theoretic problems.
While it is easy to translate a bound for $T_{H}(\varepsilon)$ into a bound for $T^*_{H}(\varepsilon)$ (in particular, establishing that $T^*_{H}(\varepsilon) \geq \poly(\varepsilon)$ for all bipartite $H$),
it is not clear if one can analogously transform a bound for $R_{E}(\varepsilon)$ into a bound for $R^*_{E}(\varepsilon)$. The first reason is that while we can average over all subsets of vertices of graphs, we can only average over ``structured'' subsets of $[n]$. This makes is hard to establish a black-box reduction/transformation between $R_{E}(\varepsilon)$ and $R^*_{E}(\varepsilon)$. The second complication is that, as we mentioned above, we do not know which equations are sparse. This makes it hard to directly relate these two quantities.
Following \cite{GHIM}, we say that $E$ is {\em abundant} if $R^*_{E}(\varepsilon) \geq \varepsilon^C$
for some $C=C(E)$. Clearly, if $E$ is abundant then it is also sparse. Gir\~{a}o, Hurley, Illingworth and Michel \cite{GHIM} asked if the converse also holds, that is,
if one can transform a polynomial bound for $R_{E}(\varepsilon)$ into a polynomial bound for $R^*_{E}(\varepsilon)$.
Our aim in this note is to prove the following.

\begin{theorem}\label{theok4}
If an invariant equation $E$ in $4$ variables is sparse, then it is also abundant. More precisely, if $R_{E}(\varepsilon) \leq \varepsilon^{-C}$ then
$R^*_{E}(\varepsilon) \geq \frac12\varepsilon^{8C}$ for all small enough $\varepsilon$.
\end{theorem}

Given the above discussion it is natural to extend the problem raised in \cite{GHIM} to all equations $E$.

\begin{problem}\label{problem2}
Is it true that for every invariant equation $E$ there is $c=c(E)$, so that for all small enough $\varepsilon$
$$
R^*_E(\varepsilon) \geq 1/R_E(\varepsilon^c)\;.
$$
\end{problem}

It is interesting to note that Varnavides \cite{V} (implicitly) gave a positive answer to Problem \ref{problem2} when $E$ is the equation $a+b=2c$.
In fact, essentially the same argument gives a positive answer to this problem for all $E$ in $3$ variables.
Hence, Problem \ref{problem2} can be considered as a generalization of Varnavides's Theorem.
Problem \ref{problem2} was also implicitly studied previously in \cite{Bloom,K}. In particular, Kosciuszko \cite{K}, extending earlier work of Schoen and Sisask \cite{SS}, gave direct lower bounds for $R^*_E$
which, thanks to \cite{KM}, are quasi-polynomially related to those of $R_E$.

The proof of Theorem \ref{theok4} is given in the next section. For the sake of completeness, and as a preparation for the proof of Theorem \ref{theok4}, we start the next section with a proof that Problem \ref{problem2} holds for equations in $3$ variables. We should point that
a somewhat unusual aspect of the proof of Theorem \ref{theok4} is that it uses a Behrend-type \cite{Behrend} geometric argument in order to
find solutions, rather than avoid them.

%\newpage

\section{Proofs}\label{seclemma}

In the first subsection below we give a concise proof of Varnavides's Theorem, namely, of the fact
that Problem \ref{problem2} has a positive answer for equations with $3$ variables.
In the second subsection we prove Theorem \ref{theok4}.

\subsection{Proof of Varnavides's theorem}
Note that for every equation $E$, there is a constant $C$ such that for every prime $p\geq Cn$ every
solution of $E$ with integers $x_i \in [n]$ over $\mathbb{F}_p$ is also a solution over $\mathbb{R}$. Since we can always find a prime $Cn \leq p \leq 2Cn$, this means that we can
assume that $n$ itself is prime\footnote{The factor $2C$ loss in the density of $S$ can be absorbed by the factor $c$ in Problem \ref{problem2}.} and count solutions over $\mathbb{F}_n$. So let $S$ be a subset of $\mathbb{F}_n$ of size $\varepsilon n$ and let $R=R_E(\varepsilon/2)$.
For $b=(b_0,b_1) \in (\mathbb{F}_n)^2$ and $x \in [R]$ let $f_{b}(x)=b_1x+b_0$ and\footnote{Since $f([R])$ is a subset of $[R]$ (rather than $S$), it might have been more accurate
to denote $f([R])$ by $f^{-1}([R])$ but we drop the $-1$ to make the notation simpler.} $f_{b}([R])=\{x \in [R]: f_{b}(x)  \in S\}$. Pick $b_0$ and $b_1$ uniformly
at random from $\mathbb{F}_n$ and note that for any $x \in [R]$ the integer $f_{b}(x)$ is uniformly distributed in $\mathbb{F}_n$. Hence,
$$
\mathbb{E}|f_{b}([R])|=\varepsilon R\;.
$$
It is also easy to see that for every $x \neq y$ the random variables $f_b(x)$ and $f_b(y)$ are pairwise independent. Hence
$$
\mathrm{Var}|f_{b}([R])| \leq \varepsilon R\;.
$$
Therefore, by Chebyshev's Inequality we have
$$
\mathbb{P}\left[|f_{b}([R])| \leq \frac{\varepsilon}{2} R\right] \leq \frac{\varepsilon R}{\varepsilon^2 R^2/4} \leq 1/2\;.
$$

In other words, at least $n^2/2$ choices of $b$ are such that $|f_{b}([R])| \geq \frac{\varepsilon}{2}R$. By our choice of $R$ this means that $f_{b}([R])$ contains $3$ distinct integers $x_1,x_2,x_3$ which satisfy $E$ and such that $f_{b}(x_i) \in S$. Note that if $x_1,x_2,x_3$ satisfy $E$ then so do $f_{b}(x_1),f_{b}(x_2),f_{b}(x_3)$. Let us denote the triple $(f_{b}(x_1),f_{b}(x_2),f_{b}(x_3))$ by $s_{b}$. We have thus
obtained $n^2/2$ solutions $s_b$ of $E$ in $S$. To conclude the proof we just need to estimate the number of times we have double counted each solution $s_{b}$. Observe that for every choice of $s_{b}=\{s_1,s_2,s_3\}$ and {\em distinct $x_1,x_2,x_3 \in [R]$}, there is exactly one choice of $b=(b_0,b_1) \in (\mathbb{F}_n)^2$ for which $b_1x_i+b_0=s_i$ for every $1 \leq i \leq 3$. Since $[R]$ contains at most $R^2$ solutions of $E$ this means that for every solution $s_1,s_2,s_3 \in S$ there are at most $R^2$ choices of $b$ for which $s_{b}=\{s_1,s_2,s_3\}$. We conclude that $S$ contains at least $n^2/2R^2$ distinct solutions, thus completing the proof.

\subsection{Proof of Theorem \ref{theok4}}

As in the proof above, we assume that $n$ is a prime and count the number of solutions of the equation $E:~\sum^4_{i=1}a_ix_i=0$ over $\mathbb{F}_n$. Let $S$ be a subset of $\mathbb{F}_n$ of size $\varepsilon n$, and let $d$ and $t$ be integers to be chosen later and let $X$ be some subset of $[t]^d$ to be chosen later as well. For every $b=(b_0,\ldots,b_d)\in (\mathbb{F}_n)^{d+1}$ and $x=(x_1,\ldots,x_d) \in X$ we use $f_b(x)$ to denote $b_0+\sum^d_{i=1}b_ix_i$ and $f_b(X)=\{x \in X: f_b(x)  \in S\}$. We call $b$ {\em good} if $|f_b(X)| \geq \varepsilon|X|/2$. We claim that at least half of all possible choices of $b$ are good. To see this,
pick $b=(b_0,\ldots,b_{d})$ uniformly at random from $(\mathbb{F}_n)^{d+1}$, and note that for any $x \in X$ the integer $f_b(x)$ is uniformly distributed in $\mathbb{F}_n$. Hence,
$$
\mathbb{E}|f_b(X)|=\varepsilon|X|\;.
$$
It is also easy to see that for every $x \neq y \in X$ the random variables
$f_b(x)$ and $f_b(y)$ are pairwise independent. Hence
$$
\mathrm{Var}|f_b(X)| \leq \varepsilon |X|\;.
$$
Therefore, by Chebyshev's Inequality we have\footnote{We will make sure $|X| \geq 8/\varepsilon$. }
\begin{equation}\label{eq1}
\mathbb{P}\left[\left|f_b(X)\right| \leq \frac{\varepsilon}{2}\left|X\right|\right] \leq 4/\varepsilon|X| \leq 1/2\;,
\end{equation}
implying that at least half of the $b$'s are good. To finish the proof we need to make sure that every such choice of a good $b$ will ``define'' a solution $s_b$ in $S$ in a way that
$s_b$ will not be identical to too many other $s_{b'}$. This will be achieved by a careful choice of $d$, $t$ and $X$.

We first choose $X$ to be the largest subset of $[t]^d$ containing no three points on one line.
We claim that
\begin{equation}\label{sizeX}
|X| \geq t^{d-2}/d\;.
\end{equation}
Indeed, for an integer $r$ let $B_r$ be the points
$x \in [t]^d$ satisfying $\sum^d_{i=1}x^2_i=r$. Then every point of $[t]^d$ lies on one such $B_r$, where $1 \leq r \leq dt^2$.
Hence, at least one such $B_r$ contains at least $t^{d-2}/d$ of the points of $[t]^d$. Furthermore, since each set $B_r$ is a subset of
a sphere, it does not contain three points on one line.

We now turn to choose $t$ and $d$. Let $C$ be such that $R_E(\varepsilon) \leq (1/\varepsilon)^C$.
Set $a=\sum^4_{i=1}|a_i|$ and pick $t$ and $d$ satisfying
\begin{equation}\label{eq2}
(1/\varepsilon)^{2C} \geq t^d \geq \left(\frac{2dt^2a^d}{\varepsilon}\right)^C\;. %\geq R_E\left(\frac{\varepsilon}{2dt^2a^d}\right) \;.
\end{equation}
Taking $t=2^{\sqrt{\log 1/\varepsilon}}$ and $d=2C\sqrt{\log 1/\varepsilon}$ satisfies\footnote{Recalling (\ref{sizeX}), we see that since $C\geq 1$ (indeed, a standard probabilistic deletion method argument shows that if an equation has $k$ variables, then $C(E) \geq 1+\frac{1}{k-2}$), we indeed have $|X| \geq 8/\varepsilon$ as we promised earlier.} the above for all small enough $\varepsilon$. Note that by (\ref{eq2}) and our choice of $C$ we have $R_E\left(\frac{\varepsilon}{2dt^2a^d}\right) \leq t^d$.

Let us call a collection of $4$ vectors $x^1,x^2,x^3,x^4 \in X$ {\em helpful} if they are distinct, and they satisfy $E$ in each coordinate, that is, for every $1 \leq i \leq d$ we have $\sum^4_{j=1}a_jx^j_i=0$.
We claim that for every good $r$, there are useful $x^1,x^2,x^3,x^4 \in f_r(X)$. To see this let $M$ denote the integers $1,\ldots,(at)^d$ and note that (\ref{sizeX}) along with the fact that $r$ is good implies that
\begin{equation}\label{eq3}
|f_r(X)| \geq \varepsilon|X|/2 \geq \frac{\varepsilon t^d}{2t^2d} = \frac{\varepsilon}{2dt^2a^d}\cdot|M|
\end{equation}
Now think of every $d$-tuple $x \in X$ as representing an integer $p(x) \in [M]$ written in base $at$. So we can also think of $f_r(X)$ as a subset of $[M]$ of density at least $\varepsilon/2dt^2a^d$. By (\ref{eq2}), we have
$$
M = (at)^d \geq t^d \geq R_E\left(\frac{\varepsilon}{2dt^2a^d}\right)\;,
$$
implying that there are {\em distinct} $x^1,x^2,x^3,x^4 \in f_r(X)$ for which $\sum^4_{j=1}a_j\cdot p(x^j)=0$. But note that since the entries of $x^1,x^2,x^3,x^4$ are from $[t]$
there is no carry when evaluating $\sum^4_{j=1}a_j\cdot p(x^j)$ in base $at$, implying that $x^1,x^2,x^3,x^4$ satisfy $E$ in each coordinate. Finally, the fact that $\sum_ja_j=0$ and that $x_i^1,x_i^2,x_i^3,x_i^4$ satisfy $E$ for each $1 \leq i \leq d$ allows us to deduce that
$$
\sum^4_{j=1}a_j \cdot f_b(x^j)=\sum^4_{j=1}a_j \cdot (b_0+\sum^d_{i=1}b_ix^j_i)=\sum^d_{i=1}b_i\cdot(\sum^4_{j=1}a_jx^j_i)=0\;,
$$
which means that $f_b(x^1),f_b(x^2),f_b(x^3),f_b(x^4)$ forms a solution of $E$. So for every good $b$, let $s_b$ be (some choice of) $f_b(x^1),f_b(x^2),f_b(x^3),f_b(x^4) \in S$ as defined above. We know from (\ref{eq1}) that
at least $n^{d+1}/2$ of all choices of $b$ are good, so we have thus obtained $n^{d+1}/2$ solutions $s_b$ of $E$ in $S$. To finish the proof we need to bound the number of times we have counted
the same solution in $S$, that is, the number of $b$ for which $s_b$ can equal a certain $4$-tuple in $S$ satisfying $E$.

Fix $s=\{s_1,s_2,s_3,s_4\}$ and recall that $s_{b}=s$ only if there is a helpful $4$-tuple $x^1,x^2,x^3,x^4$ (as defined just before equation (\ref{eq3})) such that $f_{b}(x^i)=s_i$. We claim that for every
helpful $4$-tuple $x^1,x^2,x^3,x^4$, there are at most $n^{d-2}$ choices of $b=(b_0,\ldots,b_d)$ for which $s_b=s$. Indeed recall that by our choice of $X$ the vectors $x^1,x^2,x^3$ are distinct and do not lie
on one line. Hence they are affine independent\footnote{That is, if we turn these three $d$-dimensional vectors into $(d+1)$-dimensional vectors, by adding a new coordinate whose value is $1$, we get three linearly independent vectors.} over $\mathbb{R}$. But since the entries of $x^i$ belong to $[t]$ and $t \leq 1/\varepsilon$ we see that for large enough $n$ the vectors $x^1,x^2,x^3$ are also affine
independent over $\mathbb{F}_n$. This means that the system of three linear equations
$$
b_0+b_1x^1_1+\ldots +b_dx^1_d=s_1
$$
$$
b_0+b_1x^2_1+\ldots +b_dx^2_d=s_2
$$
$$
b_0+b_1x^3_1+\ldots +b_dx^3_d=s_3
$$
(in $d+1$ unknowns $b_0,\ldots,b_d$ over $\mathbb{F}_n$) has only $n^{d-2}$ solutions, implying the desired bound on the number of choices of $b$. Since $|X| \leq t^d \leq (1/\varepsilon)^{2C}$ by (\ref{eq2}) we see that $X$ contains at most $(1/\varepsilon)^{8C}$ helpful $4$-tuples. Altogether this means that for every $s_1,s_2,s_3,s_4 \in S$ satisfying $E$, there are at most $(1/\varepsilon)^{8C}n^{d-2}$ choices of $b$ for which $s_b=s$. Since we have previously deduced that $S$ contains
at least $\frac12n^{d+1}$ solutions $s_b$, we get that $S$ contains at least $\frac12\varepsilon^{8C}n^3$ {\em distinct} solutions, as needed.

\bigskip

\noindent{\bf Acknowledgment:} I would like to thank Yuval Wigderson and an anonymous referee for helpful comments and suggestions.

\end{document}